\documentclass[11pt,reqno]{amsart}
\usepackage{inputenc}
\usepackage{enumerate}
\usepackage{amssymb}
\usepackage{amsmath}
\usepackage{mathrsfs}
\usepackage{amsthm}
\usepackage{tikz-cd}
\usepackage{mathpazo}
\usepackage{extarrows}
\usepackage{color}
\usepackage{setspace}
\usepackage[colorlinks,linkcolor = blue,anchorcolor = red,urlcolor = blue,citecolor = blue]{hyperref}
\usepackage{graphicx, subfig}
\usepackage[square, comma, sort&compress, numbers]{natbib}

\theoremstyle{definition}
\newtheorem{Def}{Definition}[section]

\newtheorem{Rem}[Def]{Remark}
\theoremstyle{plain}
\newtheorem{Thm}[Def]{Theorem}
\newtheorem{Lem}[Def]{Lemma}
\newtheorem{Pro}[Def]{Proposition}
\newtheorem{Cor}[Def]{Corollary}

\usepackage[left=2.6cm, right=2.6cm, top=3cm, bottom=3cm]{geometry}
\usepackage[all]{xy}
\usepackage{bm}
\usepackage{comment}

\def\IN{\mathbb N}\def\IC{\mathbb C}

\def\B{\mathcal B}\def\E{\mathcal E}\def\K{\mathcal K}\def\H{\mathcal H}\def\M{\mathcal M}

\def\prop{\textup{prop}}

\def\ox{\otimes}

\author[L.~Guo]{Liang Guo}
\address[L.~Guo]{Centre for Mathematical Sciences, Lund University, Box 118, Lund, 22100, Sweden}
\email{liang.guo@math.lth.se}

\author[Z.~Luo]{Zheng Luo}
\address[Z.~Luo]{College of Data Science, Jiaxing University, Jiaxing, Zhejiang, P. R. China}
\email{zluo@zjxu.edu.cn}

\author[Q.~Wang]{Qin Wang}
\address[Q.~Wang]{Research Center for Operator Algebras, School of Mathematical Sciences, East China Normal University, Shanghai, 200241, P. R. China.}
\email{qwang@math.ecnu.edu.cn}

\title[Comparing the $K$-theory of two versions of $\ell^p$-Roe algebras]{Comparing the $K$-theory of two versions of $\ell^p$-Roe algebras}
\date{\today}

\begin{document}

\begin{abstract}
We prove that for any locally finite metric space $X$ and $1<p<\infty$, the canonical inclusion from the $\ell^p$-Roe algebra with coefficient space $\ell^p(\mathbb N)$ into that with coefficient space $\ell^p(\mathbb N;\ell^2(\mathbb N))$ identifies the former with a full multiplier corner of the latter. Consequently, the inclusion induces an isomorphism on $K$-theory, answering a question of Wang, Xie, Yu, and Zhu. Moreover, at the endpoint $p=1$, we show that the two corresponding Roe algebras are not Morita equivalent whenever $X$ is infinite.
\end{abstract}

\maketitle

\section{Introduction}

Roe algebras were introduced by Roe in his work on index theory on noncompact manifolds \cite{Roe93}. They encode the large-scale geometry of metric spaces, as reflected in the rigidity of uniform Roe algebras \cite{BBFKVW22}. The $K$-theory of Roe algebras receives the higher indices of elliptic operators, and the coarse Baum--Connes conjecture predicts that this $K$-theory can be computed from the $K$-homology of the underlying space at large scales \cite{HR95, Yu00, WY20}. Banach space analogues of Roe algebras, in which Hilbert spaces are replaced by $\ell^p$-spaces, have attracted considerable attention in connection with rigidity, $K$-theory, and the $\ell^p$ coarse Baum--Connes conjecture \cite{CL18, CL21, Zha21, CN23, WXYZ24}.

Let $X$ be a locally finite discrete metric space, that is, every closed ball in $X$ is finite. For a Banach space $E$, we denote by $\B(E)$ the Banach algebra of bounded linear operators on $E$ and by $\K(E)$ its closed ideal of compact operators. For $1\le p<\infty$, every $T\in\B(\ell^p(X;E))$ is represented by a matrix $(T_{x,y})_{x,y\in X}$ with entries $T_{x,y}\in\B(E)$, and its propagation is
\[
\prop(T)=\sup\{d(x,y):T_{x,y}\ne0\}.
\]
The \emph{algebraic $\ell^p$-Roe algebra} $B^p_{\textup{alg}}(X;E)$ consists of all $T\in\B(\ell^p(X;E))$ with finite propagation and compact matrix entries, and the \emph{$\ell^p$-Roe algebra} $B^p(X;E)$ is its operator-norm closure in $\B(\ell^p(X;E))$.

The Banach space $\ell^p(X;E)$ plays the role of a geometric module over $X$, and the Roe algebra depends a priori on the coefficient space $E$. In the Hilbert space case $p=2$, every infinite-dimensional separable Hilbert space $E$ gives the same Roe algebra up to isomorphism, and its $K$-theory does not depend on this choice \cite{WY20}. For general $p$ there are two natural coefficient spaces, namely
\[
E_p=\ell^p(\IN)\qquad\text{and}\qquad E_p^{\H}=\ell^p(\IN;\H),\qquad\text{where }\H=\ell^2(\IN).
\]
Zhang and Zhou \cite{Zha21} studied the $L^p$ coarse Baum--Connes conjecture using $E_p$, whereas Wang, Xie, Yu, and Zhu \cite{WXYZ24} adopted the mixed-norm space $E_p^{\H}$ in their study of spaces that coarsely embed into $\ell^r$. We write
\[
B^p_{\textup{st}}(X)=B^p(X;E_p),\qquad B^p_{\textup{mix}}(X)=B^p(X;E_p^{\H})
\]
for the corresponding Roe algebras, where ``st'' and ``mix'' refer to the standard and mixed-norm coefficient spaces. For $p\ne2$, the spaces $E_p$ and $E_p^{\H}$ are not isomorphic, since $\ell^2$ embeds isometically into $E_p^{\H}$ but not into $E_p$ (see \cite[Corollary~2.1.6]{AK16}).

To compare them, fix a unit vector $h_0\in\H$ and define $v:E_p\to E_p^{\H}$ and $w:E_p^{\H}\to E_p$ by $v((a_n))=(a_nh_0)$ and $w((\xi_n))=(\langle\xi_n,h_0\rangle)$, where $(a_n)_{n\in\IN}\in E_p$, $(\xi_n)_{n\in\IN}\in E_p^{\H}$ and the inner product is linear in the first variable. Let $v_X:\ell^p(X;E_p)\to\ell^p(X;E_p^{\H})$ and $w_X:\ell^p(X;E_p^{\H})\to\ell^p(X;E_p)$ act coordinatewise by $(v_X\eta)(x)=v(\eta(x))$ and $(w_X\zeta)(x)=w(\zeta(x))$ for $x\in X$. Define $P\in\B(\ell^p(X;E_p^{\H}))$ and $\iota:B^p_{\textup{st}}(X)\to B^p_{\textup{mix}}(X)$ by
\begin{equation}\label{eq:corner_map}
P=v_Xw_X,\qquad \iota(T)=v_XTw_X,\quad\text{ for any } T\in B^p_{\textup{st}}(X).
\end{equation}
Our main result is the following.

\begin{Thm}\label{thm:main}
Let $X$ be a locally finite discrete metric space.
\begin{enumerate}[(1)]
\item For $1<p<\infty$, the operator $P$ is a full multiplier idempotent of $B^p_{\textup{mix}}(X)$, and $\iota$ identifies $B^p_{\textup{st}}(X)$ isometrically with the full corner $PB^p_{\textup{mix}}(X)P$. Consequently, the two Roe algebras are Morita equivalent, and $\iota$ induces an isomorphism
\[
\iota_*:K_*(B^p_{\textup{st}}(X))\xlongrightarrow{\cong}K_*(B^p_{\textup{mix}}(X)).
\]
\item For $p=1$, if $X$ is infinite, then $B^1_{\textup{st}}(X)$ and $B^1_{\textup{mix}}(X)$ are not Morita equivalent.
\end{enumerate}
\end{Thm}

Wang, Xie, Yu, and Zhu asked whether the canonical inclusion $\iota$ induces an isomorphism on $K$-theory \cite[Question~3.6]{WXYZ24}. Part (1) of Theorem \ref{thm:main} answers this question affirmatively for every $1<p<\infty$ and every locally finite discrete metric space $X$. Part (2) shows that the endpoint $p=1$ is different.

Theorem \ref{thm:main} also has an application to the coarse Baum--Connes conjecture. For a discrete metric space $X$ with bounded geometry, the $\ell^p$ coarse assembly map is defined using $\ell^p$-localization algebras over the Rips complexes of $X$ \cite{Yu97, Zha21}. We write $\mu_{\textup{st}}$ and $\mu_{\textup{mix}}$ for the assembly maps with coefficient spaces $E_p$ and $E_p^{\H}$, respectively. Combining Theorem \ref{thm:main} with the corresponding comparison of localization algebras \cite{WXYZ24} gives the following consequence.

\begin{Cor}\label{cor:CBC}
Let $X$ be a discrete metric space with bounded geometry and let $1<p<\infty$. Then the $\ell^p$ coarse assembly map $\mu_{\textup{st}}$ associated with $B^p_{\textup{st}}(X)$ is equivalent to the corresponding map $\mu_{\textup{mix}}$ associated with $B^p_{\textup{mix}}(X)$. In particular, if $X$ admits a coarse embedding into $\ell^q$ for some $1\le q<\infty$, then the $\ell^p$ coarse Baum--Connes conjecture holds for $B^p_{\textup{st}}(X)$.
\end{Cor}

The paper is organized as follows. In Section \ref{sec:full_corner}, we prove part (1) of Theorem \ref{thm:main}. Using Khintchine's inequality, we factor the finite-rank coordinate projections of $E_p^{\H}$ through $E_p$ with bounds independent of their rank. These factorizations show that $P$ is full, and the $K$-theory isomorphism then follows from Morita invariance for Banach algebras \cite{Laf02, Par09, Par15}. In Section \ref{sec:endpoint}, we prove part (2). We show that a full idempotent would force the identity of $\ell^2_n$ to factor through $\ell^1$-spaces with bounds independent of $n$, which leads to a contradiction to the lower Khintchine inequality.

\subsection*{Acknowledgment}
The first-named author is supported by the Swedish Foundation for International Cooperation in Research and Higher Education (STINT) (Grant No.~PS2025-9545). The second- and third-named authors are supported by the National Natural Science Foundation of China (NSFC) (Grant Nos.~12401155 and 12571135).

\subsection*{AI Declaration}
During the preparation of this work, the authors used Gemini 4 Argon and ChatGPT (OpenAI, GPT-6.2) as AI-assisted tools. Gemini 4 Argon pointed the authors toward the relevance of Khintchine's inequality and helped clarify the relationship between Morita equivalence and full corners for Banach algebras, particularly in the formulation and proof of Proposition \ref{pro:full_factor}. ChatGPT was used to check the logical structure of proofs, improve the exposition, and identify potential errors in the manuscript. Both tools were also used to assist with language polishing. The authors have independently reviewed and verified all mathematical arguments and take full responsibility for the content of this article.

\section{Full multiplier corners for \texorpdfstring{$1<p<\infty$}{1 < p < infty}}\label{sec:full_corner}

In this section, we prove part (1) of Theorem \ref{thm:main} using Khintchine's inequality.

For $m\in\IN$, let $\Omega_m=\{-1,1\}^m$ be equipped with the uniform probability measure. The Rademacher functions $\varepsilon_1,\ldots,\varepsilon_m : \Omega_m \to \{-1, 1\}$ are defined by
\[
\varepsilon_j(\omega)=\omega_j\quad\text{ for }\quad \omega=(\omega_1,\ldots,\omega_m)\in\Omega_m. 
\]
Khintchine's inequality states that for every $1\le s<\infty$, there exist constants $\kappa_s,c_s>0$ such that
\begin{equation}\label{eq:Khintchine}
\kappa_s\left(\sum_{j=1}^m|a_j|^2\right)^{1/2}
\le \left\|\sum_{j=1}^m a_j\varepsilon_j\right\|_{L^s(\Omega_m)}
\le c_s\left(\sum_{j=1}^m|a_j|^2\right)^{1/2}
\end{equation}
for all $m\in\IN$ and $a_1,\ldots,a_m\in\IC$.

For any $1\le p<\infty$,  $m \ge 1$ and $C > 0$, we say that the identity of $\ell^2_m$ \emph{$C$-factors through $\ell^p$} if there exist a set $I$ and bounded linear maps $U:\ell^2_m\to\ell^p(I)$ and $V:\ell^p(I)\to\ell^2_m$ such that
\[
\|VU-\textup{Id}_{\ell^2_m}\|\le \tfrac 12 \qquad\text{and}\qquad\|U\|\,\|V\|\le C.
\]
The following lemma is a standard finite-dimensional version of the fact that the Rademacher subspace is complemented in $L^p$ for $1<p<\infty$ (see for example, \cite[Proposition~6.4.2]{AK16} or \cite[Theorem~2.b.3]{LT77}), showing that the identity operators on finite-dimensional Hilbert spaces factor through finite-dimensional $\ell^p$-spaces with bounds independent of the dimension.

\begin{Lem}\label{lem:Khintchine_factor}
Let $p \in (1,\infty)$ and let $q=p/(p-1)$ be the H\"older conjugate of $p$. For every $m\in\IN$, there exist bounded linear maps
\[
U_m:\ell^2_m\longrightarrow\ell^p_{2^m},\qquad V_m:\ell^p_{2^m}\longrightarrow\ell^2_m
\]
such that
\[
V_mU_m=\textup{Id}_{\ell^2_m},\qquad \|U_m\|\le c_p,\qquad \|V_m\|\le c_q.
\]
In particular, the identity of $\ell^2_m$ $(c_pc_q)$-factors through $\ell^p$ for every $m\in\IN$.
\end{Lem}

\begin{proof}
Identify $\ell^p_{2^m}$ with $\ell^p(\Omega_m)$. Recall that the Rademacher functions $\varepsilon_1,\ldots,\varepsilon_m$ are orthonormal in $L^2(\Omega_m)$, that is,
\[
\langle\varepsilon_i,\varepsilon_j\rangle_{L^2(\Omega_m)}=2^{-m}\sum_{\omega\in\Omega_m}\varepsilon_i(\omega)\varepsilon_j(\omega)=\delta_{ij}.
\]
Using the $L^2(\Omega_m)$-inner product to recover coefficients, for $a=(a_1,\ldots,a_m)\in\ell^2_m$ and $b=(b_\omega)_{\omega\in\Omega_m}\in\ell^p(\Omega_m)$ we define
\[
U_ma=2^{-m/p}\sum_{j=1}^m a_j\varepsilon_j,\qquad (V_mb)_j=2^{m/p}\langle b,\varepsilon_j\rangle_{L^2(\Omega_m)}=2^{-m/q}\sum_{\omega\in\Omega_m}b_\omega\varepsilon_j(\omega)\quad (1\le j\le m).
\]
Then $(V_mU_ma)_j=\langle\sum_{i=1}^m a_i\varepsilon_i,\varepsilon_j\rangle_{L^2(\Omega_m)}=a_j$ for all $1\le j\le m$, and hence $V_mU_m=\textup{Id}_{\ell^2_m}$.

By Khintchine's inequality \eqref{eq:Khintchine} with $s=p$, we have
\[
\|U_ma\|_{\ell^p_{2^m}}=\left\|\sum_{j=1}^m a_j\varepsilon_j\right\|_{L^p(\Omega_m)}\le c_p\|a\|_{\ell^2_m},
\]
so $\|U_m\|\le c_p$. Moreover, by H\"older's inequality and \eqref{eq:Khintchine} with $s=q$, we obtain
\[
\begin{aligned}
\|V_mb\|_{\ell^2_m}
&=\sup_{\|a\|_{\ell^2_m}=1}\left|2^{-m/q}\sum_{\omega\in\Omega_m}b_\omega\sum_{j=1}^m\overline{a_j}\varepsilon_j(\omega)\right|\\
&\le \|b\|_{\ell^p_{2^m}}\sup_{\|a\|_{\ell^2_m}=1}\left\|\sum_{j=1}^m\overline{a_j}\varepsilon_j\right\|_{L^q(\Omega_m)}\\
&\le c_q\|b\|_{\ell^p_{2^m}}.
\end{aligned}
\]
Thus $\|V_m\|\le c_q$, completing the proof.
\end{proof}

For $k,m\in\IN$, write $\ell^p_k(\ell^2_m)=\ell^p(\{1,\ldots,k\};\ell^2_m)$ for the $\ell^p$-direct sum of $k$ copies of $\ell^2_m$. We identify $\ell^p_k(\ell^2_m)$ with the subspace of $E_p^{\H}=\ell^p(\IN;\H)$ consisting of sequences $(\xi_n)_{n\in\IN}$ such that $\xi_n\in\ell^2_m$ for $1\le n\le k$ and $\xi_n=0$ for $n>k$. Let $P_{k,m}\in\B(E_p^{\H})$ be the canonical coordinate projection onto $\ell^p_k(\ell^2_m)\subseteq E_p^{\H}$.

\begin{Pro}\label{pro:block_factor}
For every $k,m\in\IN$, $p\in (1,\infty)$, there exist finite-rank operators
\[
s_{k,m}:E_p^{\H}\longrightarrow E_p,\qquad t_{k,m}:E_p\longrightarrow E_p^{\H}
\]
such that
\begin{equation}\label{eq:skm_tkm}
t_{k,m}s_{k,m}=P_{k,m},\qquad \|s_{k,m}\|\le c_p,\qquad \|t_{k,m}\|\le c_q.
\end{equation}
Moreover, in the strong operator topology, the projections $Q_n:=P_{n,n}$ satisfy
\[
\textup{SOT-}\lim_{n\to\infty} Q_n = \textup{Id}_{E_p^{\H}},\qquad \textup{SOT-}\lim_{n\to\infty} Q_n^* = \textup{Id}_{(E_p^{\H})^*}.
\]
\end{Pro}

\begin{proof}
By Lemma \ref{lem:Khintchine_factor}, applying $U_m$ and $V_m$ to each of the $k$ Hilbert blocks gives
\[
\ell^p_k(\ell^2_m)\xrightarrow{\oplus_{j=1}^k U_m}\ell^p_{k2^m}\xrightarrow{\oplus_{j=1}^k V_m}\ell^p_k(\ell^2_m),
\]
whose composition is the identity. Let $s_{k,m}$ be the composition of the coordinate projection $E_p^{\H}\to\ell^p_k(\ell^2_m)$, the first map above, and the canonical inclusion $\ell^p_{k2^m}\hookrightarrow E_p$. Define $t_{k,m}$ similarly using the reverse maps. Then $t_{k,m}s_{k,m}=P_{k,m}$. Since the norm of an $\ell^p$-direct sum of operators equals the supremum of their norms, we obtain
\[
\|s_{k,m}\|\le c_p,\qquad \|t_{k,m}\|\le c_q.
\]

Finally, $Q_n$ converges strongly to $\textup{Id}_{E_p^{\H}}$ since the finitely supported coordinate vectors are dense in $E_p^{\H}$. Similarly, since $(E_p^{\H})^*\cong\ell^q(\IN;\H^*)$ and $1<q<\infty$, the adjoint projections $Q_n^*$ converge strongly to $\textup{Id}_{(E_p^{\H})^*}$.
\end{proof}

We next recall a consequence of Morita theory for Banach algebras developed by Lafforgue \cite{Laf02} and Paravicini \cite{Par09, Par15}. Recall that a Banach algebra $B$ is called \emph{nondegenerate} if $\overline{\textup{span}}(B B) = B$. A bounded idempotent $P \in \M(B)$ in its multiplier algebra is said to be \emph{full} if
$\overline{\textup{span}}(B P B) = B$.
By \cite[Proposition~5.3, Corollary~5.5, and Theorem~5.29]{Par09} and \cite[Th\'eor\`eme~1.2.8]{Laf02}, we have the following result.

\begin{Lem}\label{lem:corner_K}\pushQED{\qed}
Let $B$ be a nondegenerate Banach algebra and let $P \in \M(B)$ be a bounded full idempotent. Then the full corner $P B P$ is Morita equivalent to $B$, and the corner inclusion $P B P \hookrightarrow B$ induces an isomorphism
\[
K_*(P B P) \xlongrightarrow{\cong} K_*(B).\qedhere
\]
\popQED
\end{Lem}

We now construct a bounded approximate identity for $B^p_{\textup{mix}}(X)$ using the finite-rank coordinate projections $Q_n=P_{n,n}\in\B(E_p^{\H})$ from Proposition \ref{pro:block_factor}.

\begin{Lem}\label{lem:bai_Roe}
Let $X$ be a locally finite discrete metric space and let $1<p<\infty$. For each function $\lambda:X\to\IN$, define
\[
P_\lambda:=\bigoplus_{x\in X}Q_{\lambda(x)}.
\]
Then $(P_\lambda)_{\lambda\in\IN^X}$ is a contractive approximate identity for $B^p_{\textup{mix}}(X)$, where $\IN^X$ is directed by the pointwise order. Moreover, $P_\lambda\in B^p_{\textup{alg}}(X;E_p^{\H})$ and $\|P_\lambda\|\le 1$ for every $\lambda\in\IN^X$.
\end{Lem}

\begin{proof}
Each $P_\lambda$ has propagation zero and finite-rank matrix entries, with $\|P_\lambda\|\le\sup_{n\in\IN}\|Q_n\|=1$. Hence $P_\lambda\in B^p_{\textup{alg}}(X;E_p^{\H})$.

By Proposition \ref{pro:block_factor}, $Q_n$ and $Q_n^*$ converge strongly to the identity and $\sup_n\|Q_n\|=1$, the standard compactness argument gives, for every $K\in\K(E_p^{\H})$,
\[
\|(I-Q_n)K\|\to 0,\qquad \|K(I-Q_n)\|\to 0.
\]

Let $T\in B^p_{\textup{alg}}(X;E_p^{\H})$ have propagation at most $R$, and fix $\varepsilon>0$. Since $X$ is locally finite and each $T_{x,y}$ is compact, for every $x\in X$ we can choose $\lambda_0(x)\in\IN$ such that
\[
\|(I-Q_n)T_{x,y}\|<\frac{\varepsilon}{\#B(x,R)\,\#B(y,R)},\qquad \|T_{y,x}(I-Q_n)\|<\frac{\varepsilon}{\#B(x,R)\,\#B(y,R)}
\]
for all $n\ge\lambda_0(x)$ and $y\in B(x,R)$.

Now let $\lambda\ge\lambda_0$, and let $S$ denote either $T-P_\lambda T$ or $T-TP_\lambda$. Both operators have propagation at most $R$, and
\[
\|S_{x,y}\|\le\frac{\varepsilon}{\#B(x,R)\,\#B(y,R)}
\]
whenever $d(x,y)\le R$. Consequently,
\[
\sum_y\|S_{x,y}\|\le\frac{\varepsilon}{\#B(x,R)}\sum_{y\in B(x,R)}\frac{1}{\#B(y,R)}\le\varepsilon.
\]
Similarly, $\sum_x\|S_{x,y}\|\le\varepsilon$ for every $y\in X$. Applying the classical Schur test to the scalar matrix $(\|S_{x,y}\|)_{x,y\in X}$ (see \cite{Gro01}), we obtain
\[
\|T-P_\lambda T\|\le\varepsilon,\qquad \|T-TP_\lambda\|\le\varepsilon.
\]
Thus $(P_\lambda)$ is a two-sided approximate identity on $B^p_{\textup{alg}}(X;E_p^{\H})$. Since $\|P_\lambda\|\le 1$ and $B^p_{\textup{alg}}(X;E_p^{\H})$ is norm-dense in $B^p_{\textup{mix}}(X)$, the convergence extends to all of $B^p_{\textup{mix}}(X)$.
\end{proof}

We now prove part (1) of Theorem \ref{thm:main}.

\begin{proof}[Proof of Theorem \ref{thm:main} (1)]
Recall from \eqref{eq:corner_map}, since $wv=\textup{Id}_{E_p}$ and $\|v\|=\|w\|=1$, we have $P^2=P$ and $\|P\|=1$. Multiplication by $P$ preserves finite propagation and compact matrix entries, so $P$ is a multiplier idempotent of $B^p_{\textup{mix}}(X)$. For every $T\in B^p_{\textup{st}}(X)$, we have $\iota(T)=v_XTw_X\in PB^p_{\textup{mix}}(X)P$ and $w_X\iota(T)v_X=T$, whence $\|\iota(T)\|=\|T\|$. Its inverse on $PB^p_{\textup{mix}}(X)P$ is given by $S\mapsto w_XSv_X$. Thus $\iota$ identifies $B^p_{\textup{st}}(X)$ isometrically with the corner $PB^p_{\textup{mix}}(X)P$.

Write $s_n=s_{n,n}$ and $t_n=t_{n,n}$ for the finite-rank operators in Proposition \ref{pro:block_factor}, so that $t_ns_n=Q_n$, $\|s_n\|\le c_p$, and $\|t_n\|\le c_q$ for all $n\in\IN$. For each $\lambda\in\IN^X$, define
\[
C_\lambda:=\bigoplus_{x\in X}v s_{\lambda(x)},\qquad D_\lambda:=\bigoplus_{x\in X}t_{\lambda(x)}w.
\]
Both operators have propagation zero and finite-rank matrix entries, with
\[
\|C_\lambda\|\le\|v\|\sup_{n\in\IN}\|s_n\|\le c_p,\qquad \|D_\lambda\|\le\|w\|\sup_{n\in\IN}\|t_n\|\le c_q.
\]
Hence $C_\lambda,D_\lambda\in B^p_{\textup{alg}}(X;E_p^{\H})\subseteq B^p_{\textup{mix}}(X)$. Using $wv=\textup{Id}_{E_p}$ and $t_ns_n=Q_n$, we obtain
\begin{equation}\label{eq:Pnu_in_BpB}
D_\lambda P C_\lambda=\bigoplus_{x\in X}Q_{\lambda(x)}=P_\lambda.
\end{equation}
In particular, $P_\lambda\in B^p_{\textup{mix}}(X)PB^p_{\textup{mix}}(X)$. By Lemma \ref{lem:bai_Roe}, $(P_\lambda)_{\lambda\in\IN^X}$ is a bounded approximate identity for $B^p_{\textup{mix}}(X)$. Thus, for every $T\in B^p_{\textup{mix}}(X)$,
$T=\lim_\lambda P_\lambda T=\lim_\lambda D_\lambda P(C_\lambda T)$,
and since $D_\lambda,C_\lambda T\in B^p_{\textup{mix}}(X)$, we conclude that
\[
\overline{\operatorname{span}}\big(B^p_{\textup{mix}}(X)PB^p_{\textup{mix}}(X)\big)=B^p_{\textup{mix}}(X).
\]
Therefore $P$ is a full multiplier idempotent of $B^p_{\textup{mix}}(X)$. The contractive approximate identity also implies that $B^p_{\textup{mix}}(X)$ is nondegenerate, so the remaining assertions follow from Lemma \ref{lem:corner_K}.
\end{proof}

\section{The endpoint \texorpdfstring{$p=1$}{p=1}}\label{sec:endpoint}

In this section, we prove part (2) of Theorem \ref{thm:main}. Throughout this section, we assume that $X$ is an infinite locally finite discrete metric space.

In contrast to Lemma \ref{lem:Khintchine_factor}, where the identity on $\ell^2_n$ factors through $\ell^p_{2^n}$ with constants independent of $n$ for $1<p<\infty$, our argument for $p=1$ combines two opposite facts. A full idempotent forces the identity on $\ell^2_n$ to $C$-factor through $\ell^1$ with $C$ independent of $n$ (Proposition \ref{pro:full_factor}), whereas any such factorization must satisfy $C\ge c\sqrt n$ (Proposition \ref{pro:no_factor}).

For a Banach space $E$, a vector $u\in E$, and a functional $\varphi\in E^*$, we write $u\ox\varphi$ for the rank-one operator $z\mapsto\varphi(z)u$ on $E$. The following lemma identifies the minimal left ideals of the $\ell^1$-Roe algebras.

\begin{Lem}\label{lem:minimal}
If $f\in B^1_{\textup{st}}(X)$ is a nonzero idempotent with $fB^1_{\textup{st}}(X)f=\IC f$, then $f=u\ox\psi$ for some $u\in\ell^1(X;E_1)$ and $\psi\in\ell^1(X;E_1)^*$, and
\[
B^1_{\textup{st}}(X)f=\big\{z\ox\psi:z\in\ell^1(X;E_1)\big\}.
\]
The same statement holds with $B^1_{\textup{mix}}(X)$ and $E_1^{\H}$ in place of $B^1_{\textup{st}}(X)$ and $E_1$.
\end{Lem}

\begin{proof}
We only prove the statement for $B^1_{\textup{st}}(X)$, as the proof for $B^1_{\textup{mix}}(X)$ is identical. Since $f\ne 0$, finitely supported vectors are dense, and finitely supported functionals separate points, we may choose a finitely supported vector $z\in\ell^1(X;E_1)$ and a finitely supported functional $\varphi\in\ell^1(X;E_1)^*\cong\ell^\infty(X;E_1^*)$ such that $fz\ne 0$ and $\varphi\circ f\ne 0$. The operator $z\ox\varphi$ has finite propagation and finite-rank matrix entries, so $z\ox\varphi\in B^1_{\textup{st}}(X)$. Moreover,
\[
f(z\ox\varphi)f=(fz)\ox(\varphi\circ f)\in fB^1_{\textup{st}}(X)f=\IC f.
\]
Since $fz\ne 0$ and $\varphi\circ f\ne 0$, the operator $(fz)\ox(\varphi\circ f)$ is a nonzero rank-one operator. Thus $f$ has rank one and can be written as $f=u\ox\psi$ for some $u\in\ell^1(X;E_1)$ and $\psi\in\ell^1(X;E_1)^*$ with $\psi(u)=1$.

Every element of $B^1_{\textup{st}}(X)f$ therefore has the form $z\ox\psi$. Choose a finitely supported functional $\varphi_0$ such that $\varphi_0(u)=1$. Then, for every finitely supported $z$,
\[
z\ox\psi=(z\ox\varphi_0)f\in B^1_{\textup{st}}(X)f.
\]
Since finitely supported vectors are dense and $B^1_{\textup{st}}(X)f$ is closed, the result follows.
\end{proof}

\begin{Pro}\label{pro:full_factor}
Let $L$ be a Banach algebra and let $Q,R\in\M(L)$ be idempotents such that $Q$ is full in $L$, $QLQ\cong B^1_{\textup{st}}(X)$, and $RLR\cong B^1_{\textup{mix}}(X)$ as Banach algebras. Then there exists $C>0$ such that the identity of $\ell^2_n$ $C$-factors through $\ell^1$ for every $n\ge1$.
\end{Pro}

\begin{proof}
Identify $QLQ$ with $B^1_{\textup{st}}(X)$ and $RLR$ with $B^1_{\textup{mix}}(X)$ via the given isomorphisms, which changes norms only up to multiplicative constants. Since $Q$ is full in $L$, multiplying by $R$ on both sides gives
\begin{equation}\label{eq:RLR}
B^1_{\textup{mix}}(X)=\overline{\textup{span}}\big((RLQ)(QLR)\big).
\end{equation}

Fix $x_0\in X$, and choose $u\in\ell^1(X;E_1^{\H})$ and $\varphi\in\ell^1(X;E_1^{\H})^*$ supported at $x_0$ with $\varphi(u)=1$. Then $g:=u\ox\varphi\in B^1_{\textup{mix}}(X)$ is a rank-one idempotent satisfying $gLg=gB^1_{\textup{mix}}(X)g=\IC g$ and $B^1_{\textup{mix}}(X)g=\{z\ox\varphi:z\in\ell^1(X;E_1^{\H})\}$. By \eqref{eq:RLR} and $gR=Rg=g$, $g$ lies in $\overline{\textup{span}}\big((gLQ)(QLg)\big)$. Since $(gLQ)(QLg)\subseteq gLg=\IC g$ is one-dimensional, we can choose $A_0\in gLQ$ and $B_0\in QLg$ such that
\[
A_0B_0=g,\qquad gA_0=A_0,\qquad B_0g=B_0.
\]
Then $f=B_0A_0\in QLQ=B^1_{\textup{st}}(X)$ is a nonzero rank $1$ idempotent satisfying $fB^1_{\textup{st}}(X)f=\IC f$. By Lemma \ref{lem:minimal}, there exists $\psi\in\ell^1(X;E_1)^*$ such that
\[
B^1_{\textup{st}}(X)f=\big\{z\ox\psi:z\in\ell^1(X;E_1)\big\}.
\]
The maps $z\mapsto z\ox\varphi$ and $z\mapsto z\ox\psi$ identify $\ell^1(X;E_1^{\H})$ and $\ell^1(X;E_1)$ as Banach spaces with $B^1_{\textup{mix}}(X)g=RLg$ and $B^1_{\textup{st}}(X)f=QLf$, respectively, where the first identification intertwines the action of $B^1_{\textup{mix}}(X)$ with left multiplication. Moreover, since $B_0A_0=f$ and $A_0B_0=g$, right multiplication by $B_0$ and $A_0$ are mutually inverse isomorphisms between $QLf$ and $QLg$, so $QLg\cong QLf\cong \ell^1(X;E_1)$. Since left multiplication by $B\in QLR$ and $A\in RLQ$ maps $RLg$ into $QLg$ and back into $RLg$, every product $AB$ factors through $\ell^1(X;E_1)$ viewed as an operator on $\ell^1(X;E_1^{\H})$.

Identify $\ell^1(X;E_1^{\H})$ with $\ell^1(X\times\IN;\H)$, and view $\ell^2_n\subseteq\H=\ell^2(\IN)$ as the subspace spanned by the first $n$ standard basis vectors, with orthogonal projection $P_{\ell^2_n}:\H\to\ell^2_n$. Choose distinct points $x_1,x_2,\dots$ in $X$. For each $n\ge1$, define $\sigma_n:\ell^2_n\to\ell^1(X\times\IN;\H)$ and $\pi_n:\ell^1(X\times\IN;\H)\to\ell^2_n$ by
\[
(\sigma_n s)(x,i)=
\begin{cases}
s, & \text{if }(x,i)=(x_n,1),\\
0, & \text{otherwise},
\end{cases}
\qquad
\pi_n(z)=P_{\ell^2_n}\big(z(x_n,1)\big),
\]
for $s\in\ell^2_n$ and $z\in\ell^1(X\times\IN;\H)$. Then $T=\sum_{n\ge1}\sigma_n\pi_n\in B^1_{\textup{mix}}(X)$ has norm one and satisfies $\pi_nT\sigma_n=\textup{Id}_{\ell^2_n}$ for all $n\ge1$.

By \eqref{eq:RLR}, there exist $A_1,\dots,A_m\in RLQ$ and $B_1,\dots,B_m\in QLR$ such that $\|T-\sum_{j=1}^mA_jB_j\|\le 1/2$. Factoring each $A_jB_j=\Psi_j\Phi_j$ through $\ell^1(X;E_1)$ and taking the $\ell^1$-direct sum $\ell^1(I)$ of $m$ copies of $\ell^1(X;E_1)=\ell^1(X\times\IN)$, we define $U:\ell^1(X;E_1^{\H})\to\ell^1(I)$ and $V:\ell^1(I)\to\ell^1(X;E_1^{\H})$ by
\[
Uz=(\Phi_1z,\dots,\Phi_mz),\qquad V(z_1,\dots,z_m)=\sum_{j=1}^m\Psi_jz_j.
\]
Then $VU=\sum_{j=1}^mA_jB_j$, and for every $n\ge1$,
\[
\|(\pi_nV)(U\sigma_n)-\textup{Id}_{\ell^2_n}\|=\|\pi_n(VU-T)\sigma_n\|\le\frac12,\qquad\|U\sigma_n\|\,\|\pi_nV\|\le\|U\|\,\|V\|.
\]
Thus the identity of $\ell^2_n$ $C$-factors through $\ell^1$ for every $n\ge1$, with $C=\|U\|\,\|V\|$ which is independent of $n$.
\end{proof}

\begin{Rem}\label{rem:full_factor_p}
The arguments in Lemma \ref{lem:minimal} and Proposition \ref{pro:full_factor} remain valid for all $1<p<\infty$, with $\ell^1$ and $E_1$ replaced by $\ell^p$ and $E_p$. Thus, Morita equivalence requires a uniform factorization of $\textup{Id}_{\ell^2_n}$ through $\ell^p$ for all $n$. For $1<p<\infty$, this condition follows from Lemma \ref{lem:Khintchine_factor}.
\end{Rem}

The following proposition shows that Lemma~\ref{lem:Khintchine_factor} does not extend to $p=1$.

\begin{Pro}\label{pro:no_factor}
There exists a universal constant $c>0$ such that, for every $n\ge1$ and $C>0$, if the identity of $\ell^2_n$ $C$-factors through $\ell^1$, then $C\ge c\sqrt n$.
\end{Pro}

\begin{proof}
Let $I$, $U$, and $V$ be as in the definition of $C$-factorization. Replacing $V$ by $(VU)^{-1}V$, we may assume that $VU=\textup{Id}_{\ell^2_n}$ and $\|U\|\,\|V\|\le2C$. For each $i\in I$, choose $a_i,b_i\in\ell^2_n$ such that $(Ux)_i=\langle x,a_i\rangle$ for all $x\in\ell^2_n$ and $b_i$ is the image under $V$ of the $i$-th standard basis vector of $\ell^1(I)$. Then $\|b_i\|\le\|V\|$ and for any $ x\in\ell^2_n$,
\begin{equation}\label{eq:expansion}
\sum_{i\in I}|\langle x,a_i\rangle|\le\|U\|\,\|x\|,\qquad x=\sum_{i\in I}\langle x,a_i\rangle\,b_i.
\end{equation}
Let $e_1,\dots,e_n$ be the standard orthonormal basis of $\ell^2_n$. For each $\omega\in\Omega_n=\{-1,1\}^n$, the vector $x_\omega:=\sum_{k=1}^n\varepsilon_k(\omega)e_k$ has norm $\sqrt n$, and $\langle x_\omega,a\rangle=\sum_{k=1}^n\overline{a_k}\varepsilon_k(\omega)$ for every $a=(a_1,\dots,a_n)\in\ell^2_n$. By the lower bound in Khintchine's inequality \eqref{eq:Khintchine} with $s=1$, we have
\[
2^{-n}\sum_{\omega\in\Omega_n}|\langle x_\omega,a\rangle|
=\left\|\sum_{k=1}^n\overline{a_k}\varepsilon_k\right\|_{L^1(\Omega_n)}
\ge \kappa_1\|a\|_{\ell^2_n}.
\]
Applying the inequality in \eqref{eq:expansion} to $x=x_\omega$ and averaging over $\omega\in\Omega_n$ yields
\[
\kappa_1\sum_{i\in I}\|a_i\|
\le 2^{-n}\sum_{\omega\in\Omega_n}\sum_{i\in I}|\langle x_\omega,a_i\rangle|
\le \|U\|\sqrt n.
\]
Using the expansion in \eqref{eq:expansion} for $x=e_k$ and $\sup_{i\in I}\|b_i\|\le\|V\|$, we obtain
\[
\begin{aligned}
n=\sum_{k=1}^n\langle e_k,e_k\rangle
&=\sum_{k=1}^n\sum_{i\in I}\langle e_k,a_i\rangle\langle b_i,e_k\rangle
=\sum_{i\in I}\langle b_i,a_i\rangle\\
&\le\Big(\sup_{i\in I}\|b_i\|\Big)\sum_{i\in I}\|a_i\|
\le\|V\|\sum_{i\in I}\|a_i\|
\le\frac{1}{\kappa_1}\,\|U\|\,\|V\|\sqrt n
\le\frac{2C}{\kappa_1}\sqrt n,
\end{aligned}
\]
whence $C\ge(\kappa_1/2)\sqrt n$. Here the first $\leq$ above relies essentially on $p=1$.
\end{proof}

We now prove part (2) of Theorem \ref{thm:main}.

\begin{proof}[Proof of Theorem \ref{thm:main} (2)]
By Propositions \ref{pro:full_factor} and \ref{pro:no_factor}, there is no Banach algebra $L$ with idempotents $Q,R\in\M(L)$ satisfying the hypotheses of Proposition \ref{pro:full_factor}. In particular, taking $L=B^1_{\textup{mix}}(X)$ and $R=1$ shows that $B^1_{\textup{st}}(X)$ is not isomorphic to any full multiplier corner of $B^1_{\textup{mix}}(X)$.

If $B^1_{\textup{st}}(X)$ and $B^1_{\textup{mix}}(X)$ were Morita equivalent in the sense of Lafforgue \cite{Laf02} and Paravicini \cite{Par09} (see \cite[Definition~2.12]{Chu25}), then the associated linking Banach algebra
\[
L=\begin{pmatrix}B^1_{\textup{st}}(X)&E^>\\E^<&B^1_{\textup{mix}}(X)\end{pmatrix}
\]
would have diagonal matrix units $Q$ and $R=1-Q$ in $\M(L)$ satisfying $QLQ=B^1_{\textup{st}}(X)$, $RLR=B^1_{\textup{mix}}(X)$, and $\overline{\textup{span}}(LQL)=L$. This leads to a contradiction.
\end{proof}

\begin{Rem}\label{rem:K}
Theorem \ref{thm:main} (2) does not decide whether $\iota$ induces an isomorphism on $K$-theory for $p=1$, which remains open.
\end{Rem}

\bibliographystyle{alpha}
\bibliography{ref}

\end{document}